\documentclass[11pt,a4paper]{article}
\usepackage{booktabs}
\usepackage{epstopdf}
\usepackage{amsmath}
\usepackage{amsthm}
\usepackage{graphicx}
\usepackage{color}
\usepackage{amsfonts}
\usepackage{amssymb}
\usepackage{psfrag}
\usepackage{wrapfig}
\usepackage{url}
\usepackage{amscd}
\usepackage{multicol}
\usepackage{algorithm}
\usepackage{algorithmic}
\usepackage[square, comma, sort&compress, numbers]{natbib} 
\usepackage{subfig}
\usepackage{diagbox}
\usepackage{enumerate}

\newcommand{\bdes}{\begin{description}}
\newcommand{\edes}{\end{description}}
\newcommand{\bal}{\begin{align}}
\newcommand{\eal}{\end{align}}
\newcommand{\bnum}{\begin{enumerate}}
\newcommand{\enum}{\end{enumerate}}
\newcommand{\bit}{\begin{itemize}}
\newcommand{\eit}{\end{itemize}}
\newcommand{\bea}{\begin{eqnarray}}
\newcommand{\eea}{\end{eqnarray}}
\newcommand{\be}{\begin{equation}}
\newcommand{\ee}{\end{equation}}
\newcommand{\baray}{\begin{array}}
\newcommand{\earay}{\end{array}}
\newcommand{\bsry}{\begin{subarray}}
\newcommand{\esry}{\end{subarray}}
\newcommand{\bca}{\begin{cases}}
\newcommand{\eca}{\end{cases}}
\newcommand{\bcen}{\begin{center}}
\newcommand{\ecen}{\end{center}}
\newcommand{\bbm}{\begin{bmatrix}}
\newcommand{\ebm}{\end{bmatrix}}
\newcommand{\bmx}{\begin{matrix}}
\newcommand{\emx}{\end{matrix}}
\newcommand{\bpm}{\begin{pmatrix}}
\newcommand{\epm}{\end{pmatrix}}
\newcommand{\btab}{\begin{tabular}}
\newcommand{\etab}{\end{tabular}}
\theoremstyle{plain}
\newtheorem{thm}{Theorem}[section]
\newtheorem{prop}{Proposition}[section]
\newtheorem{theorem}[thm]{Theorem}

\newtheorem{lemma}{Lemma}[section]

\newtheorem{definition}{Definition}[section]
\def \mR{\mathbb{R}}

\makeatletter
\newenvironment{breakablealgorithm}
  {
   \begin{center}
     \refstepcounter{algorithm}
     \hrule height.8pt depth0pt \kern2pt
     \renewcommand{\caption}[2][\relax]{
       {\raggedright\textbf{\ALG@name~\thealgorithm} ##2\par}%
       \ifx\relax##1\relax 
         \addcontentsline{loa}{algorithm}{\protect\numberline{\thealgorithm}##2}%
       \else 
         \addcontentsline{loa}{algorithm}{\protect\numberline{\thealgorithm}##1}%
       \fi
       \kern2pt\hrule\kern2pt
     }
  }{
     \kern2pt\hrule\relax
   \end{center}
  }
\makeatother

\title{ An Adaptive High Order Method for Finding Third-Order Critical Points of Nonconvex Optimization}
\author{Xihua Zhu$^*$\quad Jiangze Han$^\dagger$\quad Bo Jiang$^\ddagger$}
\date{}

\begin{document}
\maketitle
\footnote{\hspace{-0.06in}$^*$School of Information Management and Engineering, Shanghai University of Finance and Economics, Shanghai 200433, PR China.\ Email: zhuxihua@163.sufe.edu.cn. Research of this author is supported by the GIFSUFE (Grant No. CXJJ-2019-391).}
\footnote{\hspace{-0.06in}$^\ddagger$UBC Sauder School of Business, the University of British Columbia, Vancouver, BC V6T 1Z2, Canada. Email: jiangze.han@sauder.ubc.ca.}
	\footnote{\hspace{-0.06in}$^\ddagger$Corresponding aurthor. Research Institute for Interdisciplinary Sciences, School of Information Management and Engineering, Shanghai University of Finance and Economics, Shanghai 200433, PR China. Email: isyebojiang@gmail.com. Research of this author was supported in part by NSFC Grants 11771269 and 11831002, and Program for Innovative Research Team of Shanghai University of Finance and Economics.}

\vspace{-2.0em}
\begin{abstract}
\noindent
It is well known that finding a global optimum is extremely challenging for nonconvex optimization.
There are some recent efforts \citep{anandkumar2016efficient, cartis2018second, cartis2020sharp, chen2019high} regarding the optimization methods for computing higher-order critical points, which can exclude the so-called degenerate saddle points and reach a solution with better quality. Desipte theoretical development  in \citep{anandkumar2016efficient, cartis2018second, cartis2020sharp, chen2019high}, the corresponding numerical experiments are missing. In this paper, we propose an implementable higher-order method, named adaptive high order method (AHOM), that aims to find the third-order critical points. This is achieved by solving an ``easier'' subproblem and incorporating the adaptive strategy of parameter-tuning in each iteration of the algorithm. The iteration complexity of the proposed method is established. Some preliminary numerical results are provided to show AHOM is able to escape
the degenerate saddle points, where the second-order method could possibly get stuck.

\vspace{1.0em}
\par
\vspace{1.0em}
\end{abstract}

\section{Introduction}
\begin{figure}[H]
	\begin{tabular}{p{0.03\textwidth}|cc}
		& &  \\ \hline
		\rotatebox{90}{}   &
		\subfloat
		{\includegraphics[width=0.4\textwidth, height=5.0cm]
			{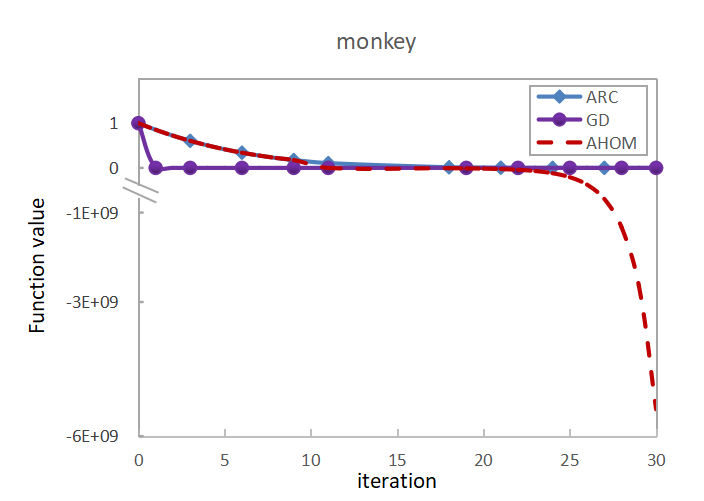}} &
		\qquad\subfloat
		{\includegraphics[width=0.4\textwidth, height=5.0cm]
			{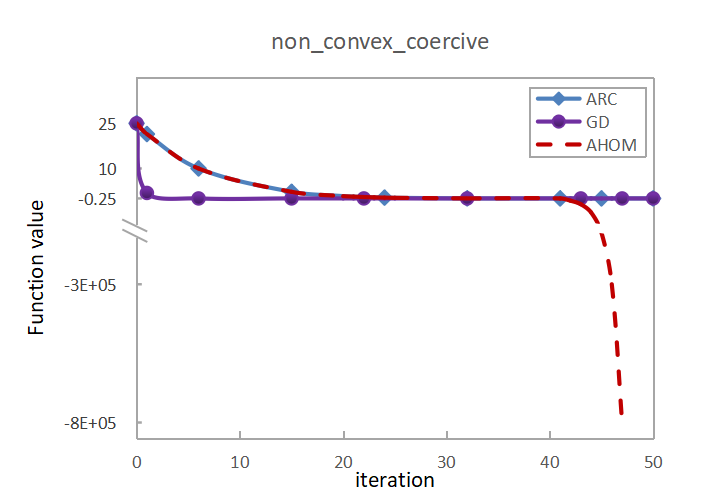}} \\
	\end{tabular}\\
	\caption{Convergence behavior (Function value v.s. Iteration) of GD, ARC and AHOM on the \textbf{monkey problem} and the \textbf{nonconvex coercive function}.}
	\label{degenerate_point}
\end{figure}

In this paper, we consider the following unconstrained optimization problem
\begin{equation}\label{Prob:main}
f^*:=\min_{x\in\mR^n}f(x),
\end{equation}
where $f$ is nonconvex and $p$-times differentiable. In recent years there have been a surge of research interest in nonconvex optimization  (see, for instance, \citep{anandkumar2016efficient, birgin2017worst, cartis2011adaptive1, cartis2011adaptive2, cartis2017improved, cartis2018second, cartis2020sharp, chen2019high, chen2019complexity, lucchi2019stochastic, gratton2019minimization, zheng2019modified, bellavia2019adaptive, gould2019convergence, curtis2018inexact, jiang2019nonconvex, martinez2017high, carmon2017lower, carmon2018accelerated, birgin2016evaluation}). However, it is well known that globally optimizing a nonconvex problem is a notoriously challenging task. Even a less challenging work of finding a local optimum
is computationally hard in the worst case \citep{nie2015hierarchy}. In fact, it is even NP-hard to check whether a critical point is a local minimizer \citep{MurtyKabadi1987}. On the other hand, the concept of critical point can be divided into
a few subcategories. Classical gradient descent type method may be stuck at a {\it first-order critical point}, i.e.,$\nabla f(x)=0$.
While algorithms incorporating second order differentiable information \citep{nesterov2006cubic} may converge to a {\it second-order critical point}, i.e., $\nabla f(x)=0$ and $\nabla^2f(x)\succeq0$, which could exclude some first-order critical points that is not local optimum. However, it is still possible that the second-order method could get stuck at the so-called {\it degenerate saddle point} (Hessian matrix has nonnegative eigenvalues with some eigenvalues equal to 0). To see this, let's consider problem \eqref{Prob:main} with two concrete objective functions:
$$
(i)\;\mbox{monkey problem: }f(x)=x_0^3-3x_0x_1^2; \quad (ii)\; \mbox{nonconvex coercive function: } f(x)=\frac{1}{3}x_0^3+\frac{1}{4}x_1^4-\frac{1}{2}x_1^2.
$$
In fact, there are two degenerate saddle points $(0,0)$ and $(0,\frac{\sqrt{3}}{3})$ in these two problems respectively. As shown in Figure \ref{degenerate_point}, the gradient descent method (GD) and the adaptive cubic regularization of Newton's method (ARC) will get stuck at these two saddle points after a few iterations by selecting $(1,0)$ and $(3,3)$ as the initial points for the two problems respectively.
To escape the degenerate saddle points, the notation of higher-order critical point was proposed in \citep{anandkumar2016efficient, cartis2018second} and the corresponding optimization algorithms were designed as well to find such higher-order critical points. By implementing these ideas to solve the monkey problem and the nonconvex coercive function, we find out in Figure \ref{degenerate_point} that by starting at the same initial points, the adaptive high order method (AHOM), which will be presented later in this paper, could indeed escape the aforementioned two degenerate saddle points, demonstrating the capability of high-order method for nonconvex optimization.

Prior to our work, there are several papers \citep{anandkumar2016efficient, cartis2018second, cartis2020sharp, chen2019high} concerning optimization methods for computing higher-order critical points. In particular, Anandkumar and Ge \citep{anandkumar2016efficient} proposed a third-order algorithm that utilizes third-order derivative information and converges to
a third-order critical point, i.e., it is a second-order point and satisfies additionally the third-order condition (see \eqref{third-order-cond} for its definition). Cartis et al. \citep{cartis2018second} presented a trust-region method with $p$-th order derivative for convexly constrained problems and it computes an $\epsilon$-approximate $q$-th ($q\ge2$) order critical points within at most $O(\epsilon^{-(q+1)})$ iterations. Later on, such iteration bound was improved to $O(\epsilon^{-\frac{q+1}{p-q+1}})$ in \citep{cartis2020sharp, chen2019high}. Despite those theoretical development, the corresponding numerical experiments are absent and the practical issue regarding the implementation of their algorithms remains to be addressed. Specifically, a nonconvex subproblem, which is NP-hard in general, needs to be globally solved in each iteration of the algorithms in \citep{cartis2018second, cartis2020sharp, chen2019high}. While Anandkumar and Ge's \citep{anandkumar2016efficient} method assumes the knowledge of problem parameters such as the Lipschitz constants of the second-order and the third-order derivatives, which is hard to estimate in practice. As a matter of fact, an algorithm that does not depend on problem parameters is often desirable in optimization. Therefore, various adaptive strategies \citep{birgin2017worst, cartis2011adaptive1, cartis2011adaptive2, cartis2020sharp, chen2019high, chen2019complexity, duchi2011adaptive, jiang2018unified, bellavia2019adaptive, zheng2019modified, lucchi2019stochastic} have been adopted to adjust the parameters in the process of iteration. In this paper, we propose an adaptive high order method (AHOM) for problem \eqref{Prob:main}, which incorporates Anandkumar and Ge's approach \citep{anandkumar2016efficient} by some adaptive strategies. In particular, the adaptation on high-order regularization term is due to the
single iteration of adaptive regularized $p$-th order method (ARp) \citep{birgin2017worst, cartis2017improved} in each step of AHOM. While the dynamic estimation on the third-order Lipschitz constant is achieved by
introducing a new successful criterion of
the third-order critical measure. It turns out our AHOM is able to solve some nonconvex machine learning problems and escape the degenerate saddle points (see Section \ref{Sec:numerical} for details).

Another merit of introducing high-order derivative information to optimization algorithms is the associated iteration complexity bounds could be improved. There are a few recent papers \citep{Baes-2009, nesterov2008accelerating, nesterov2018implementable, grapiglia2019tensorA, grapiglia2019tensorB, gasnikov2018global, jiang2018optimal, bubeck2018near} indicate that high-order derivatives indeed accelerate classical algorithms in the context of convex optimization. Similar phenomenon was also observed in the nonconvex optimization. For the unconstrained case, Nesterov \citep{nesterov2006cubic} showed that the cubic regularization of Newton method can find an $\epsilon$-approximate first-order critical point with at most $O(\epsilon^{-3/2})$ evaluations of the objective function (and its derivatives), in contrast with the evaluation complexity of $O(\epsilon^{-2})$ in the first-order method \citep{Nesterov2004Introductory}. By using up to $p$-th $(p\ge1)$ order derivatives, Birgin et al.\ \citep{birgin2017worst} first proposed ARp method, whose evaluation complexity of finding first-order critical points is improved to $O(\epsilon^{-(p+1)/p})$. Later on, Cartis et al.\ \citep{cartis2017improved} managed to adapt ARp method such that it is able to reach second-order critical points. In the mean while, the high-order method was also extended to accommodate constrained nonconvex optimization \citep{cartis2020sharp, martinez2017high, birgin2016evaluation} and non-Lipschitz nonconvex optimization \citep{chen2019high, chen2019complexity}. Since our AHOM algorithm also belongs to the category of high-order method, we show that its iteration bound improves that of the algorithm in \citep{anandkumar2016efficient}. It worths mentioning that as we perform adaptations on both the high-order regularization term and estimator of the third-order Lipschitz constant, the corresponding iteration analysis becomes more technically involved than that in \citep{anandkumar2016efficient}.

The rest of the paper is organized as follows. In Section~\ref{Sec:pre}, we introduce some preliminaries and the assumptions used throughout this paper. In Section~\ref{ARp_alg}, we propose our adaptive high order method (AHOM) for problem \eqref{Prob:main}. Section~\ref{Sec:anl} is devoted to analyzing the iteration bound of AHOM.
In Section~\ref{Sec:numerical}, we present some preliminary numerical results on solving
$\ell_2$-regularized nonconvex logistic regression problems,
where AHOM is able to escape degenerate saddle point and even occasionally converges to a point satisfying second-order sufficient condition.

\section{Preliminaries}\label{Sec:pre}
In this section, we introduce notations, various approximate critical measures and present some basic assumptions that will be used in the paper.

\subsection{Notations}
Recall that a high-order tensor is a multidimensional array. In particular first-order and second-order tensors are vectors and matrices respectively.
Throughout, we use the lower-case letters to denote vectors (e.g. $v\in\mR^n$), the capital letters to denote matrices (e.g. $M\in\mR^{n\times n}$), and the capital calligraphy letters to denote high-order tensors (e.g. $\mathcal{T}\in\mR^{n_1\times n_2\times...\times n_p}$), with subscripts of indices being their entries (e.g. $v_1,M_{i,j},\mathcal{T}_{j_1,j_2,\cdots, j_p}$). For a $p$-times differentiable function $f$, the associated $p$-th order derivative tensor is given by
\begin{equation}
\label{derivative}
\nabla^pf(x)=\left[\frac{\partial^p f(x)}{\partial x_{i_1}...\partial x_{i_p}}\right]_{i_j\in [n], \, \forall \, j},
\end{equation}
where $[n]$ denotes $\{1,...,n\}$.\par


The operations between tensor $\mathcal{T}$ and vectors $\nu^1,...,\nu^p$ yields a multi-linear form
$$
\mathcal{T}(\nu_1,...,\nu_p) = \sum_{i_1,\cdots,i_p}\mathcal{T}_{i_1\cdots i_p}\nu_{i_1}^1\cdots\nu_{i_p}^p.
$$
We say a tensor is symmetric if $\mathcal{T}_{j_1,j_2,\cdots, j_p}=\mathcal{T}_{\pi(j_1,j_2,\cdots, j_p)}$ for any permutation $\pi$ of the indices $(j_1,j_2,\cdots, j_p)$. When $\mathcal{T}$ is a symmetric tensor, we may make $\nu^1,...,\nu^p$ identical to the same vector in the above multi-linear form yielding that
$$
\mathcal{T}(\nu)^p = \sum_{i_1,\cdots,i_p}\mathcal{T}_{i_1\cdots i_p}\nu_{i_1}\cdots\nu_{i_p}.
$$
Similarly, the multi-linear form with respect to matrices $U_1\in\mR^{n\times n_1},...,U_p\in\mR^{n\times n_p}$ is defined as
$$
[\mathcal{T}(U_1,...,U_p)]_{i_1,i_2,...,i_p} = \sum_{j_1,j_2,\cdots, j_p\in[n]}\mathcal{T}_{j_1,j_2,\cdots, j_p}[U_1]_{j_1,i_1}\cdots[U_p]_{j_p,i_p},
$$
where $\mathcal{T}(U_1,...,U_p)$ itself is a $p$-th order tensor with $n_j$ being the dimension of $j$-th direction. Suppose $S$ is the projection matrix associated with subspace $\boldsymbol{S}$. We call
\begin{equation}\label{tensor-projection}
Proj_{\boldsymbol{S}}\mathcal{T} := \mathcal{T}(S,...,S)
\end{equation}
is the projection tensor of $T$ on subspace $\boldsymbol{S}$. That is,
any $p$ vectors $\nu_1,...,\nu_p$ applied to the projection tensor is equivalent to
the projections of $\nu_1,...,\nu_p$ on $\boldsymbol{S}$ applied to the original tensor:
$$
[\mathcal{T}(S,...,S)](\nu_1,...\nu_p)=\mathcal{T}(S\nu_1,...,S\nu_p).
$$


The Frobenius norm of a $p$-th order tensor $\mathcal{T}$ is: $\|\mathcal{T}\|_F=\sqrt{\sum_{j_1,j_2,\cdots, j_p\in[n]}\mathcal{T}_{j_1,j_2,\cdots, j_p}^2}$, and the spectral norm of a $p$-th order tensor is defined as
\begin{equation}
\label{spec_norm}
\|\mathcal{T}\|_{[p]}= \max_{\|\nu_1\|=...=\|\nu_p\|=1}|\mathcal{T}(\nu_1,...,\nu_p)|.
\end{equation}

For a symmetric tensor $\mathcal{T}$, the spectral norm in \eqref{spec_norm} is equivalent to $\|\mathcal{T}\|_{[p]}=\max_{\|\nu\|=1} |\mathcal{T}(\nu,...,\nu)|$. In particular, the spectral norm of a symmetric matrix $M\in\mR^{n\times n}$ is equivalent to $\|M\|_{[2]}=\max\{|\lambda_1(M)|,...,|\lambda_n(M)|\}$, where $\lambda_i(M)$ denotes the $i$-th largest eigenvalue of $M$. Note that all the matrices considered in this paper are symmetric.\par

\subsection{Lipschitz Continuous Assumption}
We assume that the $p$-th order derivative \eqref{derivative} is globally Lipschitz continuous, i.e., there exits $L_p\ge0$ such that
\begin{equation}
\label{Lipschitz}
\|\nabla^pf(x)-\nabla^pf(y)\|_{[p]}\le L_p\|x-y\|,\quad \mbox{for all}\; x,y\in R^n,
\end{equation}
where the $\|.\|_{[p]}$ is the tensor spectral norm of $p$-th order tensor given by \eqref{spec_norm}. In the rest of the paper we let $L:=\max\{\frac{L_k}{(k-1)!},~k=1,2,...,p\}$. \par
With tensor notations, the Taylor expansion of function $f(\cdot)$ at $x\in \mR^n$ can be written as:
\begin{equation}
T_p(x,s)=f(x)+\sum_{j=1}^{p}\frac{1}{j!}\nabla^jf(x)(s)^j.
\end{equation}
When $p=3$, there is a bound between $f(y)$ and its Taylor expansion (Lemma 3 in \cite{anandkumar2016efficient}):
\begin{equation}
\label{taylor_bound}
\left|f(x+s)-f(x)-\langle\nabla f(x),s\rangle -\frac{1}{2}s^{\top}\nabla^2f(x)s-\frac{1}{6}\nabla^3f(x)(s,s,s)\right|\le\frac{L_3}{24}\|s\|^4.
\end{equation}

\subsection{Approximate Critical Points}
We define the first-order and second-order critical measures of problem \eqref{Prob:main} as
\begin{equation}
\chi_{f,1}(x)\overset{def}{=}\|\nabla f(x)\|
\end{equation}
and
\begin{equation}
\chi_{f,2}(x)\overset{def}{=}\max[0,-\lambda_n(\nabla^2f(x))]
\end{equation}
respectively, where $\lambda_n(\nabla^2f(x))$ is the smallest eigenvalue of Hessian matrix $\nabla^2f(x)$. Then, a point $x$ satisfying $\chi_{f,1}(x)\le\epsilon_1$ is an $\epsilon_1$-approximate first-order critical point, and we call it an $(\epsilon_1, \epsilon_2)$-approximate second-order critical point if it further satisfies $\chi_{f,2}(x)\le\epsilon_2$.

Recall it was demonstrated in \cite{anandkumar2016efficient} that $x$ is a third-order critical point if it is a second-order critical point and
\begin{equation}\label{third-order-cond}
\nabla^3f(x)(u,u,u)=0\; \mbox{holds for any}\; u \; \mbox{that satisfies}\; u^\top\nabla^2f(x)u=0.
\end{equation}
Following the idea in \cite{anandkumar2016efficient}, we consider the eigen-subspace of the Hessian matrix below.
\begin{definition}\label{subspace}
	For any symmetric matrix $M$ with a eigen-decomposition $M=\sum_{i=1}^{n}\lambda_i v_iv_i^{\top}$, we adopt $\boldsymbol{S}_{\tau}(M)$ to denote the span of eigenvectors with eigenvalue at most $\tau$. That is
	$$\boldsymbol{S_{\tau}}(M)=span\{v_i|\lambda_i\le\tau\}.$$
\end{definition}

Now we are able to define the third-order critical measure of the objective function.

\begin{definition}[$(\beta,\kappa)$-competitive subspace and third-order critical measure]\label{comsubspace}
Given any $\beta>0$ and $\kappa>0$, let $(\beta,\kappa)$-competitive subspace $\boldsymbol{S}^x$ at point $x$ be the largest eigen-subspace $\boldsymbol{S}_{\tau}(\nabla^2f(x))$ such that $\tau\le\frac{\chi_{f,3}(x)^2}{12\kappa\beta^2}$, where
\begin{equation}
\chi_{f,3}(x)=\|Proj_{\boldsymbol{S}^x}\nabla^3f(x)\|_F
\end{equation}
is the norm of the third-order derivatives projected in this subspace. We call $\chi_{f,3}(x)$ is the third-order critical measure of $f$.
\end{definition}
Note that our notation above is slightly different from that proposed by Anandkumar and Ge \cite{anandkumar2016efficient}, where the adaptive estimator $\kappa$ is fixed as $L_3$. In contrast, we consider the $(\tau,\kappa)$-competitive subspace and third-order critical measure
to exclude the dependence on the the third-order Lipschitz parameter $L_3$.
In fact, the reason to let $\chi_{f,3}(x)=\|Proj_{\boldsymbol{S}^x}\nabla^3f(x)\|_F$ as a third-order critical measure is that condition \eqref{third-order-cond} is implied by
$\|Proj_{\boldsymbol{S}^x}\nabla^3f(x)\|_F = 0$. To see this, suppose $\|Proj_{\boldsymbol{S}^x}\nabla^3f(x)\|_F = 0$.
We observe that
$$span \{u \,| \, u^\top\nabla^2f(x)u=0 \} = \boldsymbol{S}_0(\nabla^2f(x)) \subseteq \boldsymbol{S}_\tau(\nabla^2f(x))\; \mbox{ for any} \; \tau>0$$ according to
definition \ref{subspace}. Then, for any $u\in\boldsymbol{S}_0(\nabla^2f(x))\subseteq\boldsymbol{S}_\tau(\nabla^2f(x))$, we have $u\in\boldsymbol{S}^x$. That is ${S}^xu=u$, where $S^x$ is the projection matrix associated with the subspace $\boldsymbol{S}^x$. Combining this fact with $\|Proj_{\boldsymbol{S}^x}\nabla^3f(x)\|_F=0$, we conclude that
for any $u \in \boldsymbol{S}_0(\nabla^2f(x)) \subseteq \boldsymbol{S}^x$
$$
\nabla^3f(x)(u,u,u)=\nabla^3f(x)({S}^xu,{S}^xu,{S}^xu) = 0,
$$
which is exactly the condition \eqref{third-order-cond}.

Therefore, we define the $(\epsilon_1,\epsilon_2,\epsilon_3)$-approximate critical point as follows.
\begin{definition}\label{critical-point}
We call $x$ an $(\epsilon_1,\epsilon_2,\epsilon_3)$-approximate critical point of problem \eqref{Prob:main}  if it satisfies
$$
(i)\, \chi_{f,2}(x)\le\epsilon_1,\;(ii)\, \chi_{f,2}(x)\le\epsilon_1, (iii)\, \chi_{f,3}(x)\le\epsilon_3.
$$	
\end{definition}

We end this section by presenting an algorithm named ACCS that can find a  $(\tau,\kappa)$-competitive subspace efficiently.

\begin{breakablealgorithm}
    \caption{ACCS (Algorithm for computing the $(\tau,\kappa)$-competitive subspace)}
    \label{CS}
    \begin{algorithmic}
        \STATE{$\boldsymbol{Input:}$ Hessian matrix $M=\nabla^2f(z)$, third order derivative $\mathcal{T}=\nabla^3f(z)$, approximation ratio $\beta$, adaptive parameter $\kappa$.}
        \STATE{$\boldsymbol{Output:}$ Competitive subspace $\boldsymbol{S}$ and $\chi_{f,3}(z)$}
        \STATE{\qquad Perform the eigen-decomposition of $M=\sum_{i=1}^n\lambda_iv_iv_i^{\top}$. [$\lambda_i$ is $i$-th largest eigenvalue of $M$] }
        \STATE{\qquad \textbf{for} $i=1$ to $n$ \textbf{do}\\
               \qquad\qquad Let $\boldsymbol{S}=span\{v_i,v_{i+1},...,v_n\}$. \\
               \qquad\qquad Let $\chi_{f,3}(z)=\|Proj_{\boldsymbol{S}}\mathcal{T}\|_F$\\
               \qquad\qquad \textbf{if} $\frac{\chi_{f,3}(z)^2}{12\kappa\beta^2}\ge\lambda_i$ \textbf{then}\\
               \qquad\qquad\qquad \textbf{terminate and return:} $\boldsymbol{S}$ and $\chi_{f,3}(z)$.\\
               \qquad\qquad \textbf{end if}\\
               \qquad \textbf{end for}\\
               \qquad \textbf{return} $\boldsymbol{S}=\emptyset, \chi_{f,3}(z)=0$.}
    \end{algorithmic}
\end{breakablealgorithm}

\section{Adaptive High Order Method (AHOM)}\label{ARp_alg}
This section aims to design an adaptive high order method (AHOM)
that can find a third-order critical point.
\subsection{The Single Iteration of Adaptive Regularized $p$-th Order Method (SARp) }
Before introducing the AHOM algorithm, we first present a subroutine in Algorithm \ref{ARp} that will be invoked in every iteration of AHOM. In particular, we call this subroutine SARp algorithm, which is just
a single iteration of adaptive regularized $p$-th order method (ARp) in \cite{cartis2017improved}, and requires an approximate minimization of
$$m(x_k,s,\sigma_k) :=T_p(x_k,s)+\frac{\sigma_k}{p+1}\|s\|^{p+1},
$$
where $\sigma_k$ is adaptive coefficient of the $(p+1)$-th order regularization term.
\begin{breakablealgorithm}
    \caption{Single iteration of ARp (SARp)}
    \label{ARp}
    \begin{algorithmic}
        \STATE {$\boldsymbol{Input:}$ Objective function $f$, last iterate $x_k$, regularization parameter $\sigma_k$.}
        \STATE {$\boldsymbol{Output:}$ Generated point $z_k$ and next regularization parameter $\sigma_{k+1}$.}
        \STATE {\qquad\textbf{Step 0: Initialization}. Give the constants $\theta,\eta_1,\eta_2,\gamma_1,\gamma_2,\gamma_3,\sigma_{min}$ are also given and satisfy \par\qquad\qquad\qquad\qquad\  $\theta>0, \sigma_{min}\in(0,\sigma_0], 0<\eta_1\le\eta_2<1, 0<\gamma_1<1<\gamma_2<\gamma_3.$ \par
        \qquad Compute $f(x_k)$.}
        \STATE {\qquad\textbf{Step 1: Step calculation}. Compute the step $s_k$ by approximately minimizing the model\\ \qquad $m(x_k,s,\sigma_k)$ with respect to the $s$ satisfying the following conditions
        \begin{eqnarray*}
            m(x_k,s_k,\sigma_k)&<&m(x_k,0,\sigma_k)\\
            \chi_{m,i}(x_k,s_k,\sigma_k)&\le&\theta\|s_k\|^{(p+1-i)}, \quad   (i=1,2).
        \end{eqnarray*}}
        \STATE {\qquad\textbf{Step 2: Acceptance of the trial point}. Compute $f(x_k+s_k)$ and define}
        \begin{center}
            $\rho_k=\frac{f(x_k)-f(x_k+s_k)}{T_p(x_k,0)-T_p(x_k,s_k)}$
        \end{center}
        \STATE {\qquad If $\rho_k\ge\eta_1$, then let $z_{k}=x_k+s_k$; otherwise $z_{k}=x_k$.}
        \STATE {\qquad\textbf{Step 3: Regularization parameter update}. Set
        \begin{center}
        \vspace{-4.0mm}
            \begin{displaymath}
                \sigma_{k+1}\in
                \begin{cases}
                    [\max\{\sigma_{min},\gamma_1\sigma_k \},\sigma_k] & \textrm{if $\rho\ge\eta_2$,}\\
                    [\sigma_k,\gamma_2\sigma_k] & \textrm{if $\rho_k \in [\eta_1,\eta_2)$,}\\
                    [\gamma_2\sigma_k,\gamma_3\sigma_k] & \textrm{if $\rho_k<\eta_1$.}
                \end{cases}
            \end{displaymath}
        \end{center}}
        \STATE {\qquad \textbf{Return} point $z_k$ and regularization parameter $\sigma_{k+1}$.}
    \end{algorithmic}
\end{breakablealgorithm}

We remark that the conditions in Step 1 of SARp are easily achievable by applying some existing algorithms like ARC method \cite{cartis2011adaptive1,cartis2011adaptive2}. As SARp is a single step of ARp, many useful properties of ARp can be carried over to SARp, which are summarized in the following lemma.

\begin{lemma}[\cite{cartis2017improved}, Lemma 3.1, Lemma 3.3, Lemma 3.4]
\label{stepsize} Given $x_k$,
the mechanism of SARp guarantees the following properties of the approximate minimizer $s_k$ of $m(x_k,s,\sigma_k)$.
\begin{enumerate}[\qquad\qquad\qquad\qquad\qquad\qquad(i)]
\item $T_p(x_k,0)-T_p(x_k,s_k)\ge\frac{\sigma_k}{p+1}\|s_k\|^{p+1}$,
\item $\|s_k\|\ge\left(\frac{\chi_{f,1}(x_k+s_k)}{L+\theta+\sigma_k}\right)^{\frac{1}{p}}$,
\item $\|s_k\|\ge\left(\frac{\chi_{f,2}(x_k+s_k)}{(p-1)L+\theta+p\sigma_k}\right)^{\frac{1}{p-1}}$.
\end{enumerate}
\end{lemma}
With the lemma above, we are able to prove some bounds for the critical measures and the sufficient decrease on the objective function in terms of the distance between $z_k$ and $x_k$.
\begin{prop}
\label{ARp_bound}
Suppose that $(z_k,\sigma_{k+1})=\textbf{SARp}(f,x_k,\sigma_k)$, then for all successful SARp ($\rho_k\ge\eta_1$), there have
\begin{enumerate}[\qquad\qquad\qquad\qquad\qquad\qquad(i)]
\item $\chi_{f,1}(z_k)\le(L+\theta+\sigma_k)\|z_k-x_k\|^p,$
\item $\chi_{f,2}(z_k)\le \left( (p-1)L+\theta+p\sigma_k\right)\|z_k-x_k\|^{p-1}$
\item $f(z_k)\le f(x_k)-\frac{\eta_1\sigma_{min}}{p+1}\|z_k-x_k\|^{p+1},$
\end{enumerate}
where $\eta_1$ and $\sigma_{min}$ are defined in SARp algorithm.
\end{prop}

\begin{proof} If SARp is successful, we have $z_k=x_k+s_k$,
and (i) and (ii) are just reformulations of (ii) and (iii) in Lemma \ref{stepsize}. To prove (iii), we note that
$$\rho_k=\frac{f(x_k)-f(z_k)}{T_p(x_k,0)-T_p(x_k,s_k)}\ge\eta_1$$
in successful SARp, which combined with (i) in Lemma \ref{stepsize} yields that
\begin{equation*}
\label{iteratebound}
f(x_k)-f(z_k)\ge\eta_1(T_p(x_k,0)-T_p(x_k,s_k))\ge\frac{\eta_1\sigma_{min}}{p+1}\|z_k-x_k\|^{p+1}.
\end{equation*}
\end{proof}

\subsection{The AHOM Algorithm}\label{Sec:alg}
Now we are ready to present our AHOM algorithm in Algorithm
\ref{ARp_Third}.

%

\begin{breakablealgorithm}
    \caption{Adaptive High Order Method (AHOM)}
    \label{ARp_Third}
    \begin{algorithmic}
        \STATE {$\boldsymbol{Input:}$ An initial point $x_0$, objective function $f$, accuracy levels $\epsilon_1, \epsilon_2$ ~and~$\epsilon_3$ of critical measures.}
        \STATE {$\boldsymbol{Output:}$ Solution $x_{\epsilon}$ that satisfies third-order critical measure.}
         \STATE {\qquad\textbf{Initialization}. Set regularization parameters $\sigma_0>0$, $\kappa_0>0$, and constants $0<\xi_1<1$, \\
         \qquad	$\zeta>1$, $\beta>0$. }
        \STATE {\qquad \textbf{for} $k=0,1,2,...$}
        \STATE {\qquad\qquad\textbf{Step 1: Step calculation}.\\
        \qquad\qquad\quad 1. Find $z_k$ that can decrease the objective function value by computing $$(z_k,\sigma_{k+1})=\boldsymbol{SARp}(f,x_k,\sigma_k).$$\\
        \qquad\qquad\quad 2. Search for competitive subspace $\boldsymbol{S}^{z_k}$ and third order tensor norm $\chi_{f,3}(z_k)$ by\\  \qquad\qquad\qquad computing $(\boldsymbol{S}^{z_k},\chi_{f,3}(z_k))=\boldsymbol{ACCS}(\nabla^2f(z_k),\nabla^3f(z_k),\beta,\kappa_k)$. \\
        \qquad\qquad\quad 3. {\bf Test for termination.} Evaluate $\chi_{f,i}(z_{k})$, \\
        \qquad\qquad\qquad\, {\bf if} $\chi_{f,i}(z_{k})\le\epsilon_i$, for~$i=1,2,3$, {\bf terminate} with a solution $ x_\epsilon = x_{k+1}$  \\
        \qquad\qquad\quad 4. \textbf{if} $\chi_{f,3}(z_k)\ge \beta(24\cdot\chi_{f,1}(z_k)\cdot \kappa_k^2)^{1/3}$, go to Step 2,\\
        \qquad\qquad\quad\quad \textbf{else} let $x_{k+1}=z_k$ and go to Step 3.}
        \STATE {\qquad\qquad\textbf{Step 2: Acceptance of the trial point}. \\
        \qquad\qquad Compute $u=\boldsymbol{ ATN}(\nabla^3 f(z_k),\boldsymbol{S}^{z_k},\beta)$ such that $\nabla^3f(z_k)(u,u,u)\ge \frac{\chi_{f,3}(z_k)}{\beta}$, where {\bf ATN} \\
        \qquad\qquad  is described in Algorithm \ref{ATN}. Let $\triangle_k=\frac{\chi_{f,3}(z_k)^4}{24\beta^4\kappa_k^3}$, compute $f(z_k-\frac{\chi_{f,3}(z_k)}{\beta\kappa_k}u)$ and define}
        \begin{equation}
        \label{third_success}
            \Phi_k=\frac{f(z_k)-f(z_k-\frac{\chi_{f,3}(z_k)}{\beta\kappa_k}u)}{\triangle_k}.
        \end{equation}
        \STATE {\qquad\qquad {\bf If} $\Phi_k\ge\xi_1$, let $x_{k+1}=z_k-\epsilon_k u$; {\bf otherwise} let $x_{k+1}=z_k$.}
        \STATE {\qquad\qquad\textbf{Step 3: Regularization parameter update}. Set
        \begin{center}
        \vspace{-4.0mm}
            \begin{displaymath}
                \kappa_{k+1}=
                \begin{cases}
                    \quad\zeta\kappa_k, & \textrm{if $\chi_{f,3}(z_k)\ge \beta(24\cdot\chi_{f,1}(z_k)\cdot \kappa_k^2)^{1/3}$  and $\Phi_k<\xi_1$, }\\
                    \quad\kappa_k, & \textrm{otherwise.}\\
                \end{cases}
            \end{displaymath}
        \end{center}}

            \STATE {\qquad \textbf{End for}}
    \end{algorithmic}
\end{breakablealgorithm}

AHOM algorithm utilizes the first-order, second-order and third-order derivatives to make progress, and it stops when all the three critical measures are sufficiently small, i.e., $\chi_{f,i}(z_{k})\le\epsilon_i$, for some given $\epsilon_i$ with $i=1,2,3$. The decrease of the first two order critical measures is achieved by iteratively performing SARp in Step $1$. When the third-order critical measure $\chi_{f,3}$ on the trial point $z_{k}$ is large, a descent direction $u$ will be constructed. Then a nontrivial update will be performed if the sufficient relative decrease on the objective (i.e., $\Phi_k \ge \xi_1$) further occurs. On the other hand, a step resulting in an insufficient relative decrease will be rejected by the algorithm, in the meanwhile the adaptive estimator $\kappa$ will be increased by a factor of $\zeta$. It is also possible that the third-order critical measure $\chi_{f,3}$ is already below the given tolerance but either $\chi_{f,1}$ or $\chi_{f,2}$ is still large. In this case, Step $2$ will be skipped and $x_{k+1} = z_{k}$. Furthermore, if $z_{k}$ is obtained by an unsuccessful SARp, $x_{k+1}$ actually equals to $x_{k}$ (that is $x_{k+1}$ is not updated). However, the cubic regularizer $\sigma_{k+1}$ is updated in this case, which will lead to a possible update on the next trial point $z_{k+1}$. Finally, we would like to mention that algorithm ATN in Step $2$ was proposed in \cite{anandkumar2016efficient} and is convergent by at most $2$ iterations in expectation (Theorem 7 in \cite{anandkumar2016efficient}). We present the details of ATN in the appendix for the reference of interested readers.

\section{Iteration Complexity Analysis of AHOM}\label{Sec:anl}
To provide the iteration bound for AHOM, like in \cite{birgin2017worst, cartis2017improved} we first want to define some "successful" iterations.
Since there are two regularization parameters: $\sigma$ and $\kappa$, whether they are updated successfully defines two types of "successful" iterations accordingly. We first consider the cubic regularizer $\sigma$ in SARp.
\begin{definition}
	\label{def_SARP}
	We say an iteration in AHOM is  "successful SARp" if the SARp called in Step $1$ of this iteration is successful (i.e., $\rho_j\ge\eta_1$), otherwise it is an "unsuccessful SARp" iteration. Suppose $T$ is the total number of iterations in AHOM,
    we denote by
	$$(i)~\mathcal{S}_{SARp}=\left\{0\le j\le T - 1\;|\;\rho_j\ge\eta_1\right\}$$
	the index set of all iterations such that the associated trial point $z_k$ is successful in SARp, and the complementary set including all the "unsuccessful SARp"
	iterations is denoted as
	$$(ii)~\mathcal{U}_{SARp}=\left\{0\le j\le T -1 \;|\;\rho_j<\eta_1\right\}.$$
\end{definition}

Recall that in Lemma 3.5 of \cite{cartis2017improved}, the total number of iterations in the ARp for second-order critical points can be bounded by a function of the number of successful SARp (i.e., $|\mathcal{S}_{SARp}|$). At first glance, we don't expect such bound holds true for AHOM as the iterate could possibly be updated at Step 2 of AHOM after preforming SARp, resulting a whole different sequence in contrast with that of ARp. However, we note that the universal bound in Lemma 3.2 of \cite{cartis2017improved} for the cubic regularizer $\sigma$ is still valid for $\sigma_{k+1}$ in SARp. Therefore, the same relationship between the two iteration numbers in Lemma 3.5 of \cite{cartis2017improved} is carried over to AHOM by a similar proof.


\begin{lemma}
\label{cartis_total_bound}
The mechanism of AHOM and its subroutine SARp guarantees that
\begin{equation}
T\le|\mathcal{S}_{SARp}|\left(1+\frac{|\log{\gamma_1}|}{\log{\gamma_2}}\right)+\frac{1}{\log{\gamma_2}}\log\left(\frac{\sigma_{max}}{\sigma_0}\right)
\end{equation}
where $\sigma_{\max}=\max\left\{\sigma_0,\frac{\gamma_3L(p+1)}{p(1-\eta_2)}\right\}$.
\end{lemma}
Next we consider the successful iteration defined by the update of $\kappa_{k+1}$.
\begin{definition}
	\label{def_STO}
    If an iteration in AHOM
    performs a nontrivial update $x_{j+1}=z_j-\epsilon_ju$ in Step 2, we call it an "third-order successful" iteration. Suppose $T$ is the total number of iterations in AHOM, we denote by
	$$(i)~\mathcal{S}_{third}~=\left\{0\le j\le T -1 \;|\;\Phi_j\ge\xi_1 ~and~ \chi_{f,3}(z_j)\ge \beta(24\chi_{f,1}(z_j)\kappa_j^2)^{1/3}\right\},$$
	the index set of all "third-order successful" iterations. While all the
	"third-order unsuccessful" iterations are categorized into two sets:
	$$(ii)~\mathcal{U}_{third1}=\left\{0\le j\le T - 1\;|\;\Phi_j<\xi_1~and~\chi_{f,3}(z_j)\ge\beta(24\chi_{f,1}(z_j)\kappa_j^2)^{1/3}\right\},$$
	and
	$$(iii)~\mathcal{U}_{third2}=\left\{0\le j\le T - 1\;|\;\chi_{f,3}(z_j)<\beta(24\chi_{f,1}(z_j)\kappa_j^2)^{1/3}\right\}.$$
	respectively, due to the violation on the relative decrease $\Phi$ or the third-order critical measure $\chi_{f,3}(z_j)$.
\end{definition}

According to Lemma \ref{cartis_total_bound}, it suffices to bound $|\mathcal{S}_{SARp}|$ to establish the overall iteration complexity of AHOM.
From Definition \ref{def_SARP} and \ref{def_STO}, we have the following identity: $$T -1=|\mathcal{S}_{SARp}|+|\mathcal{U}_{SARp}|=|\mathcal{S}_{third}|+|\mathcal{U}_{third1}|+|\mathcal{U}_{third2}|.$$
Consequently,
\begin{equation}
\label{SARp_bound}
\begin{aligned}
|\mathcal{S}_{SARp}|&=|\mathcal{S}_{SARp}\cap(\mathcal{S}_{third}\cup\mathcal{U}_{third1}\cup\mathcal{U}_{third2})|\\
                    &=|\mathcal{S}_{SARp}\cap\mathcal{S}_{third}|+|\mathcal{S}_{SARp}\cap\mathcal{U}_{third1}|+|\mathcal{S}_{SARp}\cap\mathcal{U}_{third2}|\\
                    &\le|\mathcal{S}_{SARp}\cap\mathcal{S}_{third}|+|\mathcal{U}_{third1}|+|\mathcal{S}_{SARp}\cap\mathcal{U}_{third2}|.
\end{aligned}
\end{equation}

Next, we shall first bound the term: $|\mathcal{U}_{third1}|$.
Before doing so, we first provide the benefit of using the third-order information.
\begin{lemma}
\label{Third_itebound}
For iteration k of AHOM, suppose $\kappa_k\ge\frac{L_3}{2-\xi_1}$, $\chi_{f,3}(z_k)\ge\beta(24\|\nabla f(z_k)\| \kappa_k^2)^{1/3}$,  $u$ is a unit vector in $\boldsymbol{S}^{z_k}$ such that $[\nabla^3f(z_k)](u,u,u)\ge \chi_{f,3}(z_k)/\beta$. Let $x_{k+1}=z_k-\frac{\chi_{f,3}(z_k)}{\beta\kappa_k}u$, then we have
\begin{equation}
f(x_{k+1})\le f(z_k)-\xi_1\frac{\chi_{f,3}(z_k)^4}{24\beta^4\kappa_k^3},
\end{equation}
i.e., iteration k is third-order successful.
\end{lemma}

\begin{proof}
Let $\epsilon=\frac{\chi_{f,3}(z_k)}{\beta\kappa_k}$, $\delta_1=\chi_{f,1}(z_k)=\|\nabla f(z_k)\|$ and $\delta_2=\max_{y \in {\bf S}^{z_k}} y^{\top} \nabla^2 f(z_k) y$, then by using \eqref{taylor_bound} we have that
$$\begin{aligned}
f(x_{k+1})&\le f(z_k)-\epsilon\nabla f(z_k)^{\top}u+\frac{\epsilon^2}{2}u^{\top}\nabla^2 f(z_k)u-\frac{\epsilon^3}{6}[\nabla^3 f(z_k)](u,u,u)+\frac{L_3\epsilon^4}{24}\|u\|^4\\
&\le f(z_k)+\delta_1\epsilon+\frac{\delta_2\epsilon^2}{2}-\frac{\|Proj_{\boldsymbol{S}^{z_k}}\nabla^3f(z_k)\|_F\epsilon^3}{6\beta}+\frac{L_3\epsilon^4}{24}\\
&=f(z_k)+\delta_1\epsilon+\frac{\delta_2\epsilon^2}{2}-\frac{\chi_{f,3}(z_k)\epsilon^3}{6\beta}+\frac{L_3\epsilon^4}{24}
\end{aligned}.$$
From the assumption $\chi_{f,3}(z_k)\ge \beta(24\delta_1\kappa_k^2)^{\frac{1}{3}}$, one has that
$$\delta_1\epsilon\le \frac{\chi_{f,3}(z_k)^3}{24\beta^3\kappa_k^2}\cdot\epsilon=\frac{\kappa_k}{24}\cdot \left(\frac{\chi_{f,3}(z_k)}{\beta\kappa_k}\right)^3\cdot\epsilon=\frac{\kappa_k\epsilon^4}{24}.$$
Furthermore the construction of the competitive subspace implies that $\delta_2\le\frac{\chi_{f,3}(z_k)^2}{12\kappa_k \beta^2}$ and thus
$$
\frac{\delta_2\epsilon^2}{2}\le\frac{\chi_{f,3}(z_k)^2}{12\kappa_k \beta^2}\cdot\frac{\epsilon^2}{2}=\frac{\kappa_k}{12}\cdot\left(\frac{\chi_{f,3}(z_k)}{\beta\kappa_k}\right)^2\cdot\frac{\epsilon^2}{2}=\frac{\kappa_k\epsilon^4}{24}.
$$
Therefore, combining the above inequalities with the assumption $\kappa_k\ge\frac{L_3}{2-\xi_1}$, which is equivalent to $\frac{2\kappa_k-L_3}{\kappa_k}\ge\xi_1$, yields that
$$\begin{aligned}
f(x_{k+1})\le f(z_k)-(2\kappa_k-L_3)\frac{\epsilon^4}{24} = f(z_k)-\frac{2\kappa_k-L_3}{\kappa_k}\cdot\frac{\chi_{f,3}(z_k)^4}{24\, \beta^4\kappa_k^3} \le f(z_k)-\xi_1\frac{\chi_{f,3}(z_k)^4}{24\beta^4\kappa_k^3},
\end{aligned}$$
 which amounts to $$\Phi_k=\frac{f(z_k)-f(z_k-\epsilon_k u)}{\triangle_k}=\frac{f(z_k)-f(x_{k+1})}{\frac{\chi_{f,3}(z_k)^4}{24\beta^4\kappa_k^3}}\ge\xi_1$$
meaning that the iteration k is a third-order successful iteration.
\end{proof}

Then it is easy to see that $\kappa_k$ has an upper bound as shown below.

\begin{lemma}
For all iteration $k$ in AHOM, we have that
\begin{equation}
\label{kappa_max}
\kappa_k\le\kappa_{\max}\overset{def}{=}\max\left\{\kappa_0,\frac{\zeta L_3}{2-\xi_1}\right\}
\end{equation}
where $\zeta>1$ and $0<\xi_1<1$.
\end{lemma}

\begin{proof}
We note that $\kappa_k$ is increased by a factor of $\zeta$ only when $\chi_{f,3}(z_k)\ge \beta(24\cdot\chi_{f,1}(z_k)\cdot \kappa_k^2)^{1/3}$ and $\Phi_k<\xi_1$. However, we have shown in Lemma \ref{Third_itebound} that $\Phi_k \ge \xi_1$ as long as $\chi_{f,3}(z_k)\ge \beta(24\cdot\chi_{f,1}(z_k)\cdot \kappa_k^2)^{1/3}$ and $\kappa_k\ge\frac{L_3}{2-\xi_1}$. Therefore, $\kappa_k$ will not be updated once it exceeds $\frac{L_3}{2-\xi_1}$.
We introduce the factor $\zeta>1$ in $\kappa_{\max}$ to accommodate case when $\kappa_k$ is only slightly less than $\frac{L_3}{2-\xi_1}$ in its last update.
\end{proof}
As a consequence,
we are able to bound the number of type $1$ unsuccessful iterations $|\mathcal{U}_{third1}|$ in AHOM.
\begin{lemma}
\label{unsec_number}
It holds that
\begin{equation} \label{Uk_bound}
|\mathcal{U}_{third1}|\le \left\lceil \frac{\log(\frac{\kappa_{\max}}{\kappa_0})}{\log\zeta} \right\rceil.\end{equation}
\end{lemma}

\begin{proof}
The updating rule of $\kappa_k$ in AHOM gives that
\begin{equation*}
\kappa_{k+1}=\zeta\kappa_k,~~k\in\mathcal{U}_{third1}, ~and~\kappa_{k+1}=\kappa_k,~~k\in(\mathcal{S}_{third}\cup\mathcal{U}_{third2}).
\end{equation*}
Thus we deduce inductively that
$$\kappa_0\zeta^{|\mathcal{U}_{third1}|}1^{|\mathcal{S}_{third}|+|\mathcal{U}_{third2}|}\le\kappa_{\max}.$$
Therefore the conclusion follows by dividing by $\kappa_0$ and then taking $\log$ on both sides.
\end{proof}

With all the above results, we are now in position to state our main complexity result below.

\begin{theorem}
Suppose algorithm AHOM starts at $x_0$, and $f$ has global min $f^*$. Then, given $\epsilon_1>0,~\epsilon_2>0$ and $\epsilon_3>0$, Algorithm \ref{ARp_Third} needs at most
$$\left\lceil\left(2\omega(f(x_0)-f^*)\max\left\{\epsilon_1^{-\frac{p+1}{p}},\epsilon_2^{-\frac{p+1}{p-1}},\epsilon_3^{-4},\epsilon_3^{-\frac{3(p+1)}{p}}\right\}+\frac{\log(\frac{\kappa_{\max}}{\kappa_0})}{\log\zeta}\right)\left(1+\frac{|\log{\gamma_1}|}{\log{\gamma_2}}\right)+\bar \Delta \right\rceil$$
iterations in total to produce an iterate $x_{\epsilon}$ such that $\chi_{f,i}(x_{\epsilon})\le \epsilon_i$, $i=1,2,3$,
where
\begin{eqnarray}
\omega\overset{def}{=}\max\Big{\{ }\frac{p+1}{\eta_1\sigma_{min}}(L+\theta+\sigma_{\max})^{\frac{p+1}{p}},\frac{p+1}{\eta_1\sigma_{min}}\left((p-1)L+\theta+p\sigma_{\max}\right)^{\frac{p+1}{p-1}},\nonumber \\
\frac{p+1}{\eta_1\sigma_{min}}(24\beta^3\kappa_{\max}^2(L+\theta+\sigma_{\max}))^{\frac{p+1}{p}},\frac{\beta\kappa_{\max}^3}{\xi_1}\Big{ \} }, \hspace{2.8cm} \label{def:omega}
\end{eqnarray}
 $\bar \Delta=\frac{1}{\log{\gamma_2}}\log\left(\frac{\sigma_{max}}{\sigma_0}\right)$, and $\kappa_{\max}$ is given by \eqref{kappa_max}.
\end{theorem}

\begin{proof} According to Lemma \ref{cartis_total_bound} and inequality \eqref{SARp_bound}, it suffices to bound the three terms: $|\mathcal{S}_{SARp}\cap\mathcal{S}_{third}|$, $|\mathcal{U}_{third1}|$ and $|\mathcal{S}_{SARp}\cap\mathcal{U}_{third2}|$ respectively. The upper bound of $|\mathcal{U}_{third1}|$ can be found in Lemma \ref{unsec_number}. To bound the other two terms, we note that the algorithm AHOM continues is due to 	
either the first-order, the second-order or the third-order critical measure is still above the given tolerance, namely,
\begin{equation}
\label{termination_condition}
(a)~\chi_{f,1}(x_{k+1})>\epsilon_1~or~(b)~\chi_{f,2}(x_{k+1})>\epsilon_2~or~(c)~\chi_{f,3}(x_{k+1})>\epsilon_3.
\end{equation}
In the following, we shall bound $|\mathcal{S}_{SARp}\cap\mathcal{S}_{third}|$, and $|\mathcal{S}_{SARp}\cap\mathcal{U}_{third2}|$ under the above three scenarios.

We first bound the term $|\mathcal{S}_{SARp}\cap\mathcal{S}_{third}|$.
Suppose the current iteration $k\in(\mathcal{S}_{SARp}\cap\mathcal{S}_{third})$, from (iii) in Proposition \ref{ARp_bound} and Proposition \ref{Third_itebound}, we know that
\begin{equation}
\label{successful_bound}
f(x_{k})-f(x_{k+1})=f(x_k)-f(z_k)+f(z_k)-f(x_{k+1})\ge\frac{\eta_1\sigma_{min}}{p+1}\|z_k-x_k\|^{p+1}+\xi_1\frac{\chi_{f,3}(z_k)^4}{\beta\kappa_k^3}.
\end{equation}
Next, we further bound the inequality \eqref{successful_bound} from below according to the three scenarios in \eqref{termination_condition}.
\begin{enumerate}
	\item In case of condition (a) in \eqref{termination_condition} holds, we deduce from \eqref{successful_bound} and part (i) in Proposition \ref{ARp_bound} that
	\begin{equation}
	\label{first}
	f(x_{k})-f(x_{k+1})\ge\frac{\eta_1\sigma_{min}}{p+1}\|z_k-x_k\|^{p+1}>\omega_1\epsilon_1^{\frac{p+1}{p}}
	\end{equation}
	where $\omega_1\overset{def}{=}\frac{\eta_1\sigma_{min}}{p+1}\left(\frac{1}{L+\theta+\sigma_{\max}}\right)^{\frac{p+1}{p}}$.
	\item In case of condition (b) in \eqref{termination_condition} holds, we deduce from \eqref{successful_bound} and part (ii) in Proposition \ref{ARp_bound} that
	\begin{equation}
	\label{second}
	f(x_{k})-f(x_{k+1})\ge\frac{\eta_1\sigma_{min}}{p+1}\|z_k-x_k\|^{p+1}>\omega_2\epsilon_2^{\frac{p+1}{p-1}}
	\end{equation}
	where $\omega_2\overset{def}{=}\frac{\eta_1\sigma_{min}}{p+1}\left(\frac{1}{(p-1)L+\theta+p\sigma_{\max}}\right)^{\frac{p+1}{p-1}}$.
	\item In case of condition (c) in \eqref{termination_condition} holds, we deduce from \eqref{successful_bound} have that
	\begin{equation}
	\label{third}
	f(x_{k})-f(x_{k+1})\ge\xi_1\frac{\chi_{f,3}(z_k)^4}{\beta\kappa_k^3}>\omega_3\epsilon_3^4
	\end{equation}
	where $\omega_3\overset{def}{=}\frac{\xi_1}{\beta\kappa_{\max}^3}$.
\end{enumerate}
Therefore, for any iteration $k\in(\mathcal{S}_{SARp}\cap\mathcal{S}_{third})$, combining \eqref{first}, \eqref{second} and \eqref{third} guarantees that:
$$f(x_{k})-f(x_{k+1})\ge\min\{\omega_1,\omega_2,\omega_3\}\min\{\epsilon_1^{\frac{p+1}{p}},\epsilon_2^{\frac{p+1}{p-1}},\epsilon_3^4\}.$$
Recalling $f^*$ is a universal lower bound of $f$, one has that
\begin{equation}
\label{sum_bound}
\begin{aligned}
f(x_0)-f^*&\ge\sum_{k=0}^{T-1}(f(x_k)-f(x_{k+1}))\\
             &\ge\sum_{k\in(\mathcal{S}_{SARp}\cap\mathcal{S}_{third})}(f(x_k)-f(x_{k+1}))\\
             &\ge\min\{\omega_1,\omega_2,\omega_3\}\min\{\epsilon_1^{\frac{p+1}{p}},\epsilon_2^{\frac{p+1}{p-1}},\epsilon_3^4\}\cdot|\mathcal{S}_{SARp}\cap\mathcal{S}_{third}|,
\end{aligned}
\end{equation}
and concludes the desired upper bound
\begin{equation}\label{SS-ub}
|\mathcal{S}_{SARp}\cap\mathcal{S}_{third}|\le\frac{f(x_0)-f^*}{\min\{\omega_1,\omega_2,\omega_3\}}\max\{\epsilon_1^{-\frac{p+1}{p}},\epsilon_2^{-\frac{p+1}{p-1}},\epsilon_3^{-4}\}.
\end{equation}

Then, we bound the number of $|\mathcal{S}_{SARp}\cap\mathcal{U}_{third2}|$.  Suppose the current iteration $k \in \mathcal{S}_{SARp}\cap\mathcal{U}_{third2}$, from (iii) in Proposition \ref{ARp_bound} and the mechanism of AHOM , we know that
\begin{equation}
\label{successful_bound_1}
f(x_{k})-f(x_{k+1})=f(x_k)-f(z_k)\ge\frac{\eta_1\sigma_{min}}{p+1}\|z_k-x_k\|^{p+1},
\end{equation}
and
\begin{equation}
\label{unsuccessful_condition}
\chi_{f,3}(z_k)<\beta(24\chi_{f,1}(z_k)\kappa_k^2)^{1/3}.
\end{equation}
Similarly, we further bound the inequality \eqref{successful_bound_1} from below according to the three scenarios in \eqref{termination_condition}. The same argument for \eqref{first} and \eqref{second} implies that they are still valid for $k \in \mathcal{S}_{SARp}\cap\mathcal{U}_{third2}$ under scenario (a) and (b) in \eqref{termination_condition} respectively. In case of third condition (c) in \eqref{termination_condition} holds, we deduce from \eqref{successful_bound_1}, \eqref{unsuccessful_condition} and part (i) in Proposition \ref{ARp_bound} have that
\begin{equation}
\label{third_2}
\begin{aligned}
f(x_{k})-f(x_{k+1})&\ge\frac{\eta_1\sigma_{min}}{p+1}\|z_k-x_k\|^{p+1}\\
                   &\ge\frac{\eta_1\sigma_{min}}{p+1}\left(\frac{\chi_{f,1}(z_k)}{L+\theta+\sigma_{\max}}\right)^{\frac{p+1}{p}}\\
                   &>\omega_4\epsilon_3^{\frac{3(p+1)}{p}}
\end{aligned}
\end{equation}
where $\omega_4\overset{def}{=}\frac{\eta_1\sigma_{min}}{p+1}\left(\frac{1}{24\beta^3\kappa_{\max}^2(L+\theta+\sigma_{\max})}\right)^{\frac{p+1}{p}}$.
Therefore, for any iteration $k\in(\mathcal{S}_{SARp}\cap\mathcal{U}_{third2})$, combining \eqref{first}, \eqref{second}, \eqref{third_2} and the argument for \eqref{SS-ub}, we have the following desired bound:
\begin{equation}\label{SU-ub}
|\mathcal{S}_{SARp}\cap\mathcal{U}_{third2}|\le\frac{f(x_0)-f^*}{\min\{\omega_1,\omega_2,\omega_4\}}\max\{\epsilon_1^{-\frac{p+1}{p}},\epsilon_2^{-\frac{p+1}{p-1}},\epsilon_3^{-\frac{3(p+1)}{p}}\}.
\end{equation}
Therefore, combining inequalities \eqref{SARp_bound}, \eqref{Uk_bound}, \eqref{SS-ub} and \eqref{SU-ub} altogether, one has that
$$\begin{aligned}
|\mathcal{S}_{SARp}|&\le\frac{f(x_0)-f^*}{\min\{\omega_1,\omega_2,\omega_3\}}\max\left\{\epsilon_1^{-\frac{p+1}{p}},\epsilon_2^{-\frac{p+1}{p-1}},\epsilon_3^{-4}\right\}\\
                    &\quad+\frac{f(x_0)-f^*}{\min\{\omega_1,\omega_2,\omega_4\}}\max\left\{\epsilon_1^{-\frac{p+1}{p}},\epsilon_2^{-\frac{p+1}{p-1}},\epsilon_3^{-\frac{3(p+1)}{p}}\right\}
                    +\frac{\log(\frac{\kappa_{\max}}{\kappa_0})}{\log\zeta}\\
                    &\le2\frac{f(x_0)-f^*}{\min\{\omega_1,\omega_2,\omega_3,\omega_4\}}\max\left\{\epsilon_1^{-\frac{p+1}{p}},\epsilon_2^{-\frac{p+1}{p-1}},\epsilon_3^{-4},\epsilon_3^{-\frac{3(p+1)}{p}}\right\}+\frac{\log(\frac{\kappa_{\max}}{\kappa_0})}{\log\zeta}
\end{aligned}$$
Finally, by invoking Lemma \ref{cartis_total_bound}, the total number of iterations $T$ can be upper bounded by
$$\begin{aligned}
T&\le|\mathcal{S}_{SARp}|\left(1+\frac{|\log{\gamma_1}|}{\log{\gamma_2}}\right)+\frac{1}{\log{\gamma_2}}\log\left(\frac{\sigma_{max}}{\sigma_0}\right)\\
 &\le\left(2\omega(f(x_0)-f^*)\max\{\epsilon_1^{-\frac{p+1}{p}},\epsilon_2^{-\frac{p+1}{p-1}},\epsilon_3^{-4},\epsilon_3^{-\frac{3(p+1)}{p}}\}+\frac{\log(\frac{\kappa_{\max}}{\kappa_0})}{\log\zeta}\right)\left(1+\frac{|\log{\gamma_1}|}{\log{\gamma_2}}\right)\\
 &\quad+\frac{1}{\log{\gamma_2}}\log\left(\frac{\sigma_{max}}{\sigma_0}\right),
\end{aligned}$$
where $\omega$ is defined in \eqref{def:omega}.
\end{proof}


\section{Numerical Experiments}\label{Sec:numerical}
In this section, we show the performance of our algorithm for solving the following nonconvex logistic regression problem:
\begin{equation}\label{nonconvex-logistic}
\min\limits_{w\in\mR^d} \frac{1}{2}\sum_{i=1}^{n}\left(\frac{1}{1+e^{-w^{\top}x_i}}-y_i\right)^2 +\frac{\alpha}{2}\|w\|^2,
\end{equation}
where $\{(x_i,y_i)\}_{i=1}^n$ is a collection of data samples with $y_i$ labeled as
$0$ or $1$, and the regularization parameter $\alpha=10^{-5}$. In contrast with the standard logistic regression, the loss in \eqref{nonconvex-logistic} is quantified as
the square of the difference between the logistic function $\frac{1}{1+e^{-w^{\top}x_i}}$ and the observed outcome $y_i$, which thus is a nonconvex function. In fact, there has been some similar nonconvex loss function proposed and studied in \cite{ghadimi2019generalized, Mason1999Nips}. It worths mentioning that both the loss functions in \eqref{nonconvex-logistic} and
\cite{ghadimi2019generalized, Mason1999Nips} belong to a broader function class named sigmoid function,  and there is an optimization model tailored for sigmoid function called sigmoidal programming in \cite{udell2013maximizing}.

Our experiments are conducted on 6 datasets all come from LIBSVM available at\\
 https://www.csie.ntu.edu.tw/\textasciitilde cjlin/libsvmtools/datasets/binary.html. The summary of those datasets are shown in Table \ref{Tab-datasets}.
\begin{table}[H]
\caption{Statistics of datasets.} \label{Tab-datasets}\vspace{-1em}
\begin{center}
\begin{tabular}{c|c|c} \hline
Dataset             & Number of Samples & Dimension \\ \hline
\textit{a1a}        & 1,605             & 123       \\
\textit{phishing} & 11,055 & 68 \\
\textit{sonar}      & 208               & 60        \\
\textit{splice}     & 1,000             & 60        \\
\textit{svmguide1}  & 3,089             & 4         \\
\textit{svmguide3}  & 1,243             & 22        \\
\hline
\end{tabular}
\end{center}
\end{table}
\vspace{-2.0em}
We apply AHOM algorithm with $p=3$ in the experiments and the subroutine SARp in each iteration reduces to the adaptive cubic regularization of Newton's method (ARC). The parameters for the subroutine SARp are set as
$\sigma_0=2.0, \sigma_{\min}=10^{-16}, \gamma_1=0.5, \gamma_2=1.1, \gamma_3=2.0, \eta_1=0.1, \eta_2=0.9$, which are the same as that of the benchmark algorithm: the adaptive cubic regularized Newton's method (ARC).
For the parameters in the main procedure of AHOM,
we set $\xi_1=10^{-9}, \zeta=1.1$, $\kappa_0 = 10^{-6}$ and $\beta=20$. Particularly, we apply the so-called Lanczos \cite{cartis2011adaptive1} process to approximately solve the subproblem $\min_{s\in\mR^n} m(x_k,s,\sigma_k)$ in SARp.

We compare our AHOM method (p=3) with two second order methods including: ARC and the trust region method (TR). We adopt the implementation of ARC and TR in the public package$^1$ with the default parameters except that the full batch rather than the subsampled batch of the component functions is taken. In ARC, the ``initial\_penalty\_parameter" $\sigma_0=2.0$ and the ``penalty\_decrease\_multiplier" $\gamma_2=1.1$ (all parameters are same as those in SARp), the initial radius and the max radius of trust region in TR is $5$ and $10^4$ respectively.
\footnote{$^1$ https://github.com/dalab/subsampled cubic regularization.}

Recall that the first-order, second-order and third-order critical measures are given by $\chi_{f,1}(x)=\|\nabla f(x)\|$, $\chi_{f,2}(x)=\max[0,-\lambda_n(\nabla^2f(x))]$, $\chi_{f,3}(x)=\|Proj_{\boldsymbol{S}(x)}\nabla^3f(x)\|_F$ respectively. We set equal error tolerances of $10^{-6}$ for these three measures. Then an approximate third-order critical point satisfies (i) $\chi_{f,1}(x)\leq 10^{-6}$, (ii) $\chi_{f,2}(x)\leq10^{-6}$, (iii)$\chi_{f,3}(x)\le 10^{-6}$, which is also the stopping criterion for AHOM. In addition, ARC and TR stops when (i) and (ii) are satisfied.\par

\begin{figure}[htbp]
	\begin{tabular}{p{0.03\textwidth}|ccc}
		& &  \\ \hline
		\rotatebox{90}{}&
		\subfloat
		{\includegraphics[width=0.3\textwidth, height=4.5cm]
			{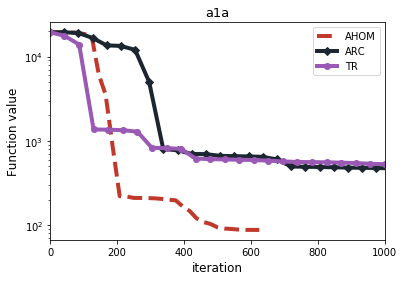}}&
		\subfloat
		{\includegraphics[width=0.3\textwidth, height=4.5cm]
			{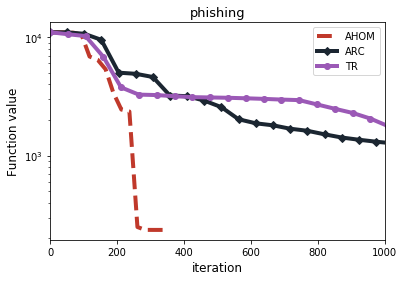}}&
		\subfloat
		{\includegraphics[width=0.3\textwidth, height=4.5cm]
			{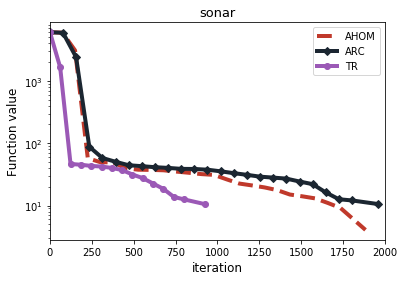}} \\
		
		\rotatebox{90}{}  &
		\subfloat
		{\includegraphics[width=0.3\textwidth, height=4.5cm]
			{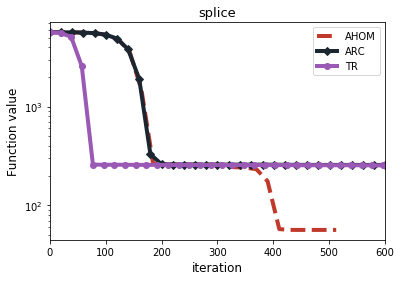}} &
		\subfloat
		{\includegraphics[width=0.3\textwidth, height=4.5cm]
			{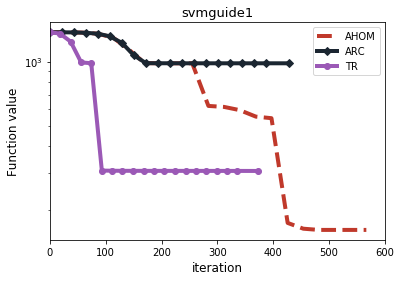}}&
		\subfloat
		{\includegraphics[width=0.3\textwidth, height=4.5cm]
			{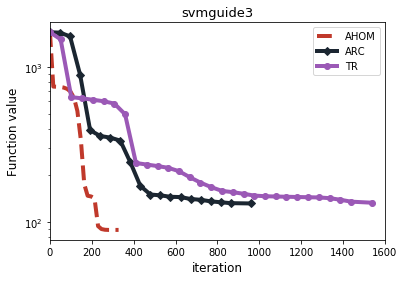}}\\
	\end{tabular}\\
	\caption{\small Performance of three algorithms on  \textbf{regularized nonconvex logistic regression} (Loss v.s. Iteration)}
	\label{iter_1e5}
\end{figure}

\begin{figure}[htbp]
	\begin{tabular}{p{0.01\textwidth}|ccc}
		& &  \\ \hline
		\rotatebox{90}{}&
		\subfloat
		{\includegraphics[width=0.3\textwidth, height=4.5cm]
			{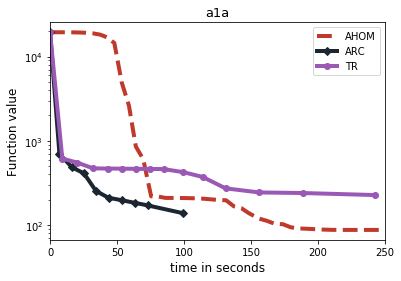}}&
		\subfloat
		{\includegraphics[width=0.3\textwidth, height=4.5cm]
			{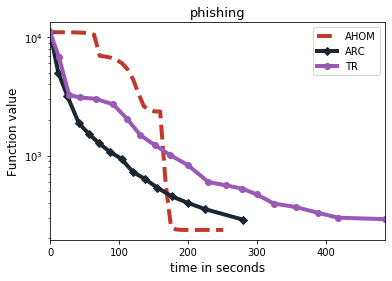}}&
		\subfloat
		{\includegraphics[width=0.3\textwidth, height=4.5cm]
			{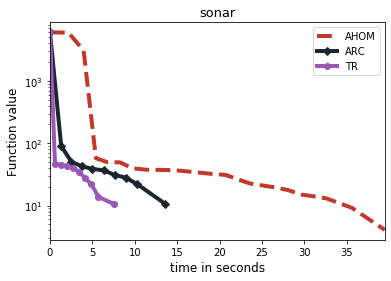}} \\
		
		\rotatebox{90}{}  &
		\subfloat
		{\includegraphics[width=0.3\textwidth, height=4.5cm]
			{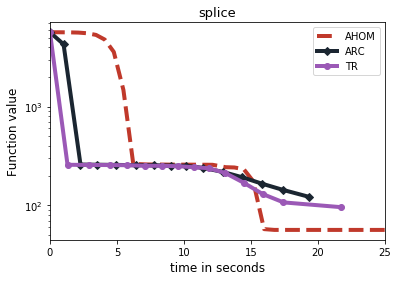}}&
		\subfloat
		{\includegraphics[width=0.3\textwidth, height=4.5cm]
			{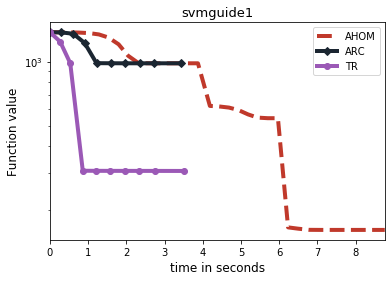}}&
		\subfloat
		{\includegraphics[width=0.3\textwidth, height=4.5cm]
			{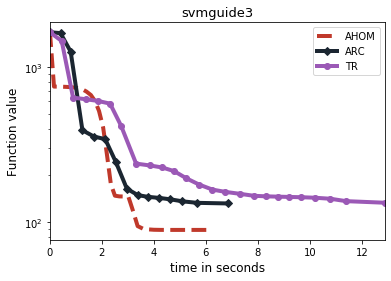}} \\
	\end{tabular}\\
	\caption{\small Performance of three algorithms on  \textbf{regularized nonconvex logistic regression } (Loss v.s. Time)}
	\label{time_1e5}
\end{figure}

We plot figures to visualize the performance of AHOM and the two benchmark methods for solving problem \eqref{nonconvex-logistic}. In particular,
function value versus iterations and function values versus time of the three algorithms are shown in Figure \ref{iter_1e5} and Figure \ref{time_1e5} respectively. We can see that AHOM is able to converge to a better solution for all the six datasets, while the other two methods may get stuck at some lower-order critical point.
This indicates that using third-order information can really help to escape some degenerate saddle points. The detailed information about function value and the three critical measures of the output points of the three algorithms is presented in Table \ref{convergence_1e5}, where we can see that AHOM even occasionally converges to a point satisfying second-order sufficient condition (see the rows for the dataset ``splice'' and ``svmguide1'').

\begin{table}[H]
	\caption{Function value and the critical measures at the converged points}\vspace{-1em}
	\label{convergence_1e5}
	\begin{center}
		\begin{tabular}{c|cccc}
			\hline
			Algorithm     & $f(x)$                              & $\|\nabla f(x)\|$
			& $\lambda_{min}(\nabla^2f(x))$       & $\|Proj_{S}\nabla^3f(x)\|_F$                \\
			\hline
			\multicolumn{5}{c}{a1a}                                                                           \\
			\hline
			\textit{AHOM} & \multicolumn{1}{c}{87.6710}         & \multicolumn{1}{c}{1.8122e-10}
			& \multicolumn{1}{c}{1e-05}           & \multicolumn{1}{c}{1.6990e-13}              \\
			\textit{ARC}  & \multicolumn{1}{c}{138.1366}        & \multicolumn{1}{c}{4.8090e-07}
			& \multicolumn{1}{c}{1e-05}           & \multicolumn{1}{c}{--}                      \\
			\textit{TR}   & \multicolumn{1}{c}{227.2495}        & \multicolumn{1}{c}{3.8878e-07}
			& \multicolumn{1}{c}{4.8523e-06}      & \multicolumn{1}{c}{--}                      \\
			\hline
			\multicolumn{5}{c}{phishing}                                                                      \\
			\hline
			\textit{AHOM} & \multicolumn{1}{c}{236.1797}        & \multicolumn{1}{c}{5.5960e-08 }
			& \multicolumn{1}{c}{1e-05}           & \multicolumn{1}{c}{0}                       \\
			\textit{ARC}  & \multicolumn{1}{c}{287.9241}        & \multicolumn{1}{c}{9.4991e-07}
			& \multicolumn{1}{c}{8.5778e-06}      & \multicolumn{1}{c}{--}                      \\
			\textit{TR}   & \multicolumn{1}{c}{291.0511}        & \multicolumn{1}{c}{6.1220e-07}
			& \multicolumn{1}{c}{1e-05}           & \multicolumn{1}{c}{--}                      \\
			\hline
			\multicolumn{5}{c}{sonar}                                                                         \\
			\hline
			\textit{AHOM} & \multicolumn{1}{c}{4.0587}          & \multicolumn{1}{c}{4.1230e-07}
			& \multicolumn{1}{c}{9.3519e-06}      & \multicolumn{1}{c}{2.5533e-15}              \\
			\textit{ARC}  & \multicolumn{1}{c}{10.5456}         & \multicolumn{1}{c}{4.9407e-07}
			& \multicolumn{1}{c}{9.8599e-06}      & \multicolumn{1}{c}{--}                      \\
			\textit{TR}   & \multicolumn{1}{c}{10.5456}         & \multicolumn{1}{c}{7.1172e-08}
			& \multicolumn{1}{c}{9.8573e-06}      & \multicolumn{1}{c}{--}                      \\
			\hline
			\multicolumn{5}{c}{splice}                                                                        \\
			\hline
			\textit{AHOM} & \multicolumn{1}{c}{56.2595}         & \multicolumn{1}{c}{3.4501e-13}
			& \multicolumn{1}{c}{0.3029}          & \multicolumn{1}{c}{0}                       \\
			\textit{ARC}  & \multicolumn{1}{c}{116.7087}        & \multicolumn{1}{c}{7.6571e-07}
			& \multicolumn{1}{c}{6.7674e-05}      & \multicolumn{1}{c}{--}                      \\
			\textit{TR}   & \multicolumn{1}{c}{95.6435}         & \multicolumn{1}{c}{2.0769e-07}
			& \multicolumn{1}{c}{7.0064e-05}      & \multicolumn{1}{c}{--}                      \\
			\hline
			\multicolumn{5}{c}{svmguide1}                                                                     \\
			\hline
			\textit{AHOM} & \multicolumn{1}{c}{161.2606}        & \multicolumn{1}{c}{8.6737e-10}
			& \multicolumn{1}{c}{4.0163}          & \multicolumn{1}{c}{0}                       \\
			\textit{ARC}  & \multicolumn{1}{c}{982.5015}        & \multicolumn{1}{c}{7.2893e-07}
			& \multicolumn{1}{c}{4.4215e-05}      & \multicolumn{1}{c}{--}                      \\
			\textit{TR}   & \multicolumn{1}{c}{305.6594}        & \multicolumn{1}{c}{2.9699e-08}
			& \multicolumn{1}{c}{3.7649e-05}      & \multicolumn{1}{c}{--}                      \\
			\hline
			\multicolumn{5}{c}{svmguide3}                                                                     \\
			\hline
			\textit{AHOM} & \multicolumn{1}{c}{89.1117}         & \multicolumn{1}{c}{7.7390e-07}
			& \multicolumn{1}{c}{1e-05}           & \multicolumn{1}{c}{0}                       \\
			\textit{ARC}  & \multicolumn{1}{c}{131.9918}        & \multicolumn{1}{c}{8.1708e-07}
			& \multicolumn{1}{c}{1e-05}           & \multicolumn{1}{c}{--}                      \\
			\textit{TR}   & \multicolumn{1}{c}{133.4102}        & \multicolumn{1}{c}{8.0748e-07}
			& \multicolumn{1}{c}{1e-05}           & \multicolumn{1}{c}{--}                      \\
			\hline
		\end{tabular}
	\end{center}
\end{table}

\section*{Acknowledgments}
We would like to thank Professor Qi Deng at Shanghai University of Finance and Ecomoics for the discussion on the numerical experiment of this paper.

\bibliographystyle{apalike}

\begin{thebibliography}{}
	
	\bibitem[Anandkumar and Ge, 2016]{anandkumar2016efficient}
	Anandkumar, A. and Ge, R. (2016).
	\newblock Efficient approaches for escaping higher order saddle points in
	non-convex optimization.
	\newblock In {\em Conference on Learning Theory}, pages 81--102.
	
	\bibitem[Baes, 2009]{Baes-2009}
	Baes, M. (2009).
	\newblock Estimate sequence methods: extensions and approximations.
	\newblock {\em Institute for Operations Research, ETH, Zürich, Switzerland}.
	
	\bibitem[Bellavia et~al., 2019]{bellavia2019adaptive}
	Bellavia, S., Gurioli, G., Morini, B., and Toint, P.~L. (2019).
	\newblock Adaptive regularization algorithms with inexact evaluations for
	nonconvex optimization.
	\newblock {\em SIAM Journal on Optimization}, 29(4):2881--2915.
	
	\bibitem[Birgin et~al., 2016]{birgin2016evaluation}
	Birgin, E.~G., Gardenghi, J., Mart{\'\i}nez, J.~M., Santos, S.~A., and Toint,
	P.~L. (2016).
	\newblock Evaluation complexity for nonlinear constrained optimization using
	unscaled {KKT} conditions and high-order models.
	\newblock {\em SIAM Journal on Optimization}, 26(2):951--967.
	
	\bibitem[Birgin et~al., 2017]{birgin2017worst}
	Birgin, E.~G., Gardenghi, J., Mart{\'\i}nez, J.~M., Santos, S.~A., and Toint,
	P.~L. (2017).
	\newblock Worst-case evaluation complexity for unconstrained nonlinear
	optimization using high-order regularized models.
	\newblock {\em Mathematical Programming}, 163(1-2):359--368.
	
	\bibitem[Bubeck et~al., 2018]{bubeck2018near}
	Bubeck, S., Jiang, Q., Lee, Y.~T., Li, Y., and Sidford, A. (2018).
	\newblock Near-optimal method for highly smooth convex optimization.
	\newblock {\em arXiv preprint arXiv:1812.08026}.
	
	\bibitem[Carmon et~al., 2017]{carmon2017lower}
	Carmon, Y., Duchi, J.~C., Hinder, O., and Sidford, A. (2017).
	\newblock Lower bounds for finding stationary points
	\uppercase\expandafter{\romannumeral1}.
	\newblock {\em Mathematical Programming}, pages 1--50.
	
	\bibitem[Carmon et~al., 2018]{carmon2018accelerated}
	Carmon, Y., Duchi, J.~C., Hinder, O., and Sidford, A. (2018).
	\newblock Accelerated methods for nonconvex optimization.
	\newblock {\em SIAM Journal on Optimization}, 28(2):1751--1772.
	
	\bibitem[Cartis et~al., 2011a]{cartis2011adaptive1}
	Cartis, C., Gould, N.~I., and Toint, P.~L. (2011a).
	\newblock Adaptive cubic regularisation methods for unconstrained optimization.
	part \uppercase\expandafter{\romannumeral1}: motivation, convergence and
	numerical results.
	\newblock {\em Mathematical Programming}, 127(2):245--295.
	
	\bibitem[Cartis et~al., 2011b]{cartis2011adaptive2}
	Cartis, C., Gould, N.~I., and Toint, P.~L. (2011b).
	\newblock Adaptive cubic regularisation methods for unconstrained optimization.
	part \uppercase\expandafter{\romannumeral2}: worst-case function-and
	derivative-evaluation complexity.
	\newblock {\em Mathematical programming}, 130(2):295--319.
	
	\bibitem[Cartis et~al., 2017]{cartis2017improved}
	Cartis, C., Gould, N.~I., and Toint, P.~L. (2017).
	\newblock Improved second-order evaluation complexity for unconstrained
	nonlinear optimization using high-order regularized models.
	\newblock {\em arXiv preprint arXiv:1708.04044}.
	
	\bibitem[Cartis et~al., 2018]{cartis2018second}
	Cartis, C., Gould, N.~I., and Toint, P.~L. (2018).
	\newblock Second-order optimality and beyond: Characterization and evaluation
	complexity in convexly constrained nonlinear optimization.
	\newblock {\em Foundations of Computational Mathematics}, 18(5):1073--1107.
	
	\bibitem[Cartis et~al., 2020]{cartis2020sharp}
	Cartis, C., Gould, N.~I., and Toint, P.~L. (2020).
	\newblock Sharp worst-case evaluation complexity bounds for arbitrary-order
	nonconvex optimization with inexpensive constraints.
	\newblock {\em SIAM Journal on Optimization}, 30(1):513--541.
	
	\bibitem[Chen and Toint, 2020]{chen2019high}
	Chen, X. and Toint, P.~L. (2020).
	\newblock High-order evaluation complexity for convexly-constrained
	optimization with {Non-Lipschitzian} group sparsity terms.
	\newblock {\em Mathematical Programming, to appear}.
	
	\bibitem[Chen et~al., 2019]{chen2019complexity}
	Chen, X., Toint, P.~L., and Wang, H. (2019).
	\newblock Complexity of partially separable convexly constrained optimization
	with {Non-Lipschitzian} singularities.
	\newblock {\em SIAM Journal on Optimization}, 29(1):874--903.
	
	\bibitem[Curtis et~al., 2018]{curtis2018inexact}
	Curtis, F.~E., Robinson, D.~P., and Samadi, M. (2018).
	\newblock An inexact regularized {Newton} framework with a worst-case iteration
	complexity of for nonconvex optimization.
	\newblock {\em IMA Journal of Numerical Analysis}, 39(3):1296--1327.
	
	\bibitem[Duchi et~al., 2011]{duchi2011adaptive}
	Duchi, J., Hazan, E., and Singer, Y. (2011).
	\newblock Adaptive subgradient methods for online learning and stochastic
	optimization.
	\newblock {\em Journal of Machine Learning Research}, 12(Jul):2121--2159.
	
	\bibitem[Gasnikov et~al., 2018]{gasnikov2018global}
	Gasnikov, A., Kovalev, D., Mohhamed, A., and Chernousova, E. (2018).
	\newblock The global rate of convergence for optimal tensor methods in smooth
	convex optimization.
	\newblock {\em arXiv preprint arXiv:1809.00382}.
	
	\bibitem[Ghadimi et~al., 2019]{ghadimi2019generalized}
	Ghadimi, S., Lan, G., and Zhang, H. (2019).
	\newblock Generalized uniformly optimal methods for nonlinear programming.
	\newblock {\em Journal of Scientific Computing}, 79(3):1854--1881.
	
	\bibitem[Gould et~al., 2019]{gould2019convergence}
	Gould, N.~I., Rees, T., and Scott, J.~A. (2019).
	\newblock Convergence and evaluation-complexity analysis of a regularized
	tensor-{Newton} method for solving nonlinear least-squares problems.
	\newblock {\em Computational Optimization and Applications}, 73(1):1--35.
	
	\bibitem[Grapiglia and Nesterov, 2019a]{grapiglia2019tensorB}
	Grapiglia, G. and Nesterov, Y. (2019a).
	\newblock Tensor methods for finding approximate stationary points of convex
	functions.
	\newblock {\em arXiv preprint arXiv:1907.07053}.
	
	\bibitem[Grapiglia and Nesterov, 2019b]{grapiglia2019tensorA}
	Grapiglia, G. and Nesterov, Y. (2019b).
	\newblock Tensor methods for minimizing functions with h{\"o}lder continuous
	higher-order derivatives.
	\newblock {\em arXiv preprint arXiv:1904.12559}.
	
	\bibitem[Gratton et~al., 2019]{gratton2019minimization}
	Gratton, S., Simon, E., and Toint, P.~L. (2019).
	\newblock Minimization of nonsmooth nonconvex functions using inexact
	evaluations and its worst-case complexity.
	\newblock {\em arXiv preprint arXiv:1902.10406}.
	
	\bibitem[Jiang et~al., 2019]{jiang2019nonconvex}
	Jiang, B., Lin, T., Ma, S., and Zhang, S. (2019).
	\newblock Structured nonconvex and nonsmooth optimization: Algorithms and
	iteration complexity analysis.
	\newblock {\em Computational Optimization and Applications}, 72:115--157.
	
	\bibitem[Jiang et~al., 2020a]{jiang2018unified}
	Jiang, B., Lin, T., and Zhang, S. (2020a).
	\newblock A unified adaptive tensor approximation scheme to accelerate
	composite convex optimization.
	\newblock {\em SIAM Journal on Optimization, to appear}.
	
	\bibitem[Jiang et~al., 2020b]{jiang2018optimal}
	Jiang, B., Wang, H., and Zhang, S. (2020b).
	\newblock An optimal high-order tensor method for convex optimization.
	\newblock {\em Mathematics of Operations Research, to appear}.
	
	\bibitem[Lucchi and Kohler, 2019]{lucchi2019stochastic}
	Lucchi, A. and Kohler, J. (2019).
	\newblock A stochastic tensor method for non-convex optimization.
	\newblock {\em arXiv preprint arXiv:1911.10367}.
	
	\bibitem[Mart{\'\i}nez, 2017]{martinez2017high}
	Mart{\'\i}nez, J.~M. (2017).
	\newblock On high-order model regularization for constrained optimization.
	\newblock {\em SIAM Journal on Optimization}, 27(4):2447--2458.
	
	\bibitem[Mason et~al., 1999]{Mason1999Nips}
	Mason, L., Baxter, J., Bartlett, P., and Frean, M. (1999).
	\newblock Boosting algorithms as gradient descent infunction space.
	\newblock In {\em Advances in neural information processing systems}, pages
	512--518.
	
	\bibitem[Murty and Kabadi, 1987]{MurtyKabadi1987}
	Murty, K. and Kabadi, S. (1987).
	\newblock Some np-complete problems in quadratic and nonlinear programming.
	\newblock {\em Mathematical Programming}, 39:117--129.
	
	\bibitem[Nesterov, 2004]{Nesterov2004Introductory}
	Nesterov, Y. (2004).
	\newblock Introductory lectures on convex optimization a basic course.
	\newblock {\em Applied Optimization}, 87(5):xviii,236.
	
	\bibitem[Nesterov, 2008]{nesterov2008accelerating}
	Nesterov, Y. (2008).
	\newblock Accelerating the cubic regularization of {Newton'}s method on convex
	problems.
	\newblock {\em Mathematical Programming}, 112(1):159--181.
	
	\bibitem[Nesterov, 2019]{nesterov2018implementable}
	Nesterov, Y. (2019).
	\newblock Implementable tensor methods in unconstrained convex optimization.
	\newblock {\em Mathematical Programming}, published
	onlone:doi:10.1007/s10107--019--01449--1.
	
	\bibitem[Nesterov and Polyak, 2006]{nesterov2006cubic}
	Nesterov, Y. and Polyak, B.~T. (2006).
	\newblock Cubic regularization of {Newton} method and its global performance.
	\newblock {\em Mathematical Programming}, 108(1):177--205.
	
	\bibitem[Nie, 2015]{nie2015hierarchy}
	Nie, J. (2015).
	\newblock The hierarchy of local minimums in polynomial optimization.
	\newblock {\em Mathematical Programming}, 151(2):555--583.
	
	\bibitem[Udell and Boyd, 2013]{udell2013maximizing}
	Udell, M. and Boyd, S. (2013).
	\newblock Maximizing a sum of sigmoids.
	\newblock {\em Optimization and Engineering}, pages 1--25.
	
	\bibitem[Zheng and Zheng, 2019]{zheng2019modified}
	Zheng, Y. and Zheng, B. (2019).
	\newblock A modified adaptive cubic regularization method for large-scale
	unconstrained optimization problem.
	\newblock {\em arXiv preprint arXiv:1904.07440}.
	
\end{thebibliography}

\appendix
\renewcommand{\appendixname}{Appendix~\Alph{section}}
\section{Algorithm for Approx Tensor Norm}
The details of algorithm ATN in Step $2$ of AHOM are shown as follows.

\begin{breakablealgorithm}
	\caption{Approximate Tensor Norms (ATN)}
	\label{ATN}
	\begin{algorithmic}
		\STATE{$\boldsymbol{Input:}$ Tensor $T$, subspace $\boldsymbol{S}$, constant $\beta$.}
		\STATE{$\boldsymbol{Output:}$ unit vector $u\in\boldsymbol{S}$ such that $T(u,u,u)\ge\|Proj_{\boldsymbol{S}}T\|_F/\beta$.}
		\STATE{\qquad \textbf{repeat}\\
			\qquad\qquad Let $u$ be a random standard Gaussian in subspace $\boldsymbol{S}$.\\
			\qquad \textbf{until} $|T(u,u,u)|\ge\|Proj_{\boldsymbol{S}}T\|_F/\beta$.\\
			\qquad \textbf{return} $u$ if $T(u,u,u)>0$ and $-u$ otherwise.}
	\end{algorithmic}
\end{breakablealgorithm}
The above algorithm was proposed in \cite{anandkumar2016efficient} and it aims to find a vector $u\in\boldsymbol{S}$ such that the value $T(u,u,u)$ is an approximation of $\|Proj_{\boldsymbol{S}}T\|_F$. The following theorem reveals that this algorithm converges by at most $2$ iterations in expectation.
\begin{theorem}[\cite{anandkumar2016efficient} Theorem 7]
There is a universal constant B such that the expected number of iterations of Algorithm ATN is at most 2, and the output of ATN is a unit vector $u$ that satisfies $T(u,u,u)\ge\|Proj_{\boldsymbol{S}}T\|_F/\beta$ for $\beta={Bn^{1.5}}$.
\end{theorem}


\end{document}